\documentclass{amsart}

\usepackage{amssymb,latexsym,amsxtra,amscd,verbatim,ifthen}
\usepackage{amscd}
\usepackage{hyperref}

\usepackage[normalem]{ulem}

\DeclareMathOperator{\Spec}{Spec}

\DeclareMathOperator{\character}{char}

\newcommand{\scrm}{{\mathfrak m}}

\newcommand{\ZZ}{{\mathbb Z}}

\newcommand{\CC}{{\mathbb C}}
\newcommand{\NN}{{\mathbb N}}
\newcommand{\PP}{{\mathbb P}}
\newcommand{\LL}{{\mathcal L}}
\newcommand{\F}{{\mathcal F}}
\newcommand{\G}{{\mathcal G}}
\newcommand{\calE}{{\mathcal E}}

\newcommand{\calH}{{\mathcal H}}
\newcommand{\calI}{{\mathcal I}}
\newcommand{\calK}{{\mathcal K}}
\newcommand{\calN}{{\mathcal N}}

\newcommand{\calP}{{\mathcal P}}

\newcommand{\OO}{{\mathcal O}}

\numberwithin{equation}{section}

\newtheorem{theorem}[equation]{Theorem}

\newtheorem{lemma}[equation]{Lemma}

\theoremstyle{definition}

\newtheorem{definition}[equation]{Definition}

\theoremstyle{remark}
\newtheorem{remark}[equation]{Remark}

\begin{document}

\title{Fujita's Conjecture and Frobenius amplitude}

\author{Dennis S. Keeler}
     \thanks{
     Partially supported by an NSA Young Investigators Grant.}  
\address{Department of Mathematics \\ Miami University \\ Oxford, OH 45056 }
\email{keelerds@miamioh.edu}
\urladdr{http://www.users.miamioh.edu/keelerds/}
 


%
%

\begin{abstract}
We prove a version of Fujita's Conjecture in arbitrary characteristic,
generalizing results of K.E.~Smith. Our methods use the Frobenius
morphism, but avoid tight closure theory. We also obtain 
versions of Fujita's Conjecture for coherent sheaves with
certain ampleness properties.
\end{abstract}

\maketitle

\section{Introduction}

Fujita's Conjecture is a deceptively simple open question in classical
algebraic geometry. Given a smooth complex projective variety $X$
of dimension $d$ and
an ample line bundle $\LL$, the conjecture predicts that
\begin{enumerate}
\item The line bundle $\omega_X \otimes \LL^{m}$ is generated
by global sections for $m \geq d + 1$.
\item The line bundle $\omega_X \otimes \LL^{m}$ is very ample
for $m \geq d + 2$.
\end{enumerate}

While the conjecture was stated in some form more than two decades
ago \cite{Fuj-orig}, thus far the global generation
conjecture has only been proven for $\dim X \leq 4$ \cite{EinLaz,Kaw}.
Also, over $\CC$, many other ``Fujita Conjecture type'' theorems
have been proven. We direct the reader to \cite[Section~10.4]{PAG}
for a partial summary.

Since these proofs rely on the Kodaira Vanishing Theorem
and its generalizations, it was surprising when K.E.~Smith proved
in arbitrary characteristic,
via tight closure theory, that if
 $\LL$ is ample and \emph{generated by global sections}, then
$\omega_X \otimes \LL^{d+1}$ is generated by global sections
\cite{Smith-FujConj1}. 
This result was recovered in \cite{Hara-FujConj},
via a characteristic $p$ analogue of multiplier ideals.
Using tight closure methods, Smith also proved that if
$\LL$ is very ample and $(X,\LL) \neq (\PP^n, \OO(1))$, 
then $\omega_X \otimes \LL^d$
is generated by global sections \cite{Smith-FujConj2}.
(Note that Smith and Hara allowed $F$-rational singularities,
but we will remain in the smooth case.)

Using characteristic $p$ methods, but staying within the realm
of algebraic geometry, we will prove

\begin{theorem}\label{thm:FujitaConj}
Let $X$ be a projective scheme of pure dimension $d$, smooth over a field $k$
of arbitrary characteristic. Let $\LL$ be
an ample, globally generated line bundle and let $\calH$ be an ample line bundle. Then
\begin{enumerate}
\item\label{fuj1} $\omega_X \otimes \LL^{d} \otimes \calH$ is generated by global sections.
\item\label{fuj2} $\omega_X \otimes \LL^{d + 1} \otimes \calH$ is very ample.
\end{enumerate}
\end{theorem}

As far as the author knows, this is the greatest generality in which
Fujita's Conjecture has been proven in characteristic $p$. Also, this
is the first characteristic $p$ version of the very ampleness part
of the Conjecture. On the other hand, the theorem above can be deduced
easily over $\CC$, using the Kodaira Vanishing Theorem and Castelnuovo-Mumford
regularity \cite[Example~1.8.23]{PAG}. While we still use 
regularity, we of course avoid Kodaira Vanishing. Thus we hope
our method will be useful in other situations where one must avoid Kodaira Vanishing.
(For those wishing to avoid vanishing theorems altogether,  one may 
prove Theorem~\ref{thm:FujitaConj} as in \cite[Remark~10.4.5]{PAG}
when $\calH = \LL$ is very ample,
regardless of $\character k$.)

In \cite{AraK}, Arapura undertakes a thorough investigation of so-called
$F$-ample coherent sheaves and generalizes the definition to
characteristic $0$ (see Definition~\ref{def:F-ample}). 
Our method easily generalizes Theorem~\ref{thm:FujitaConj}\eqref{fuj1}
to the case where $\calH$ is any $F$-ample coherent sheaf 
(see Theorem~\ref{thm:mainF-ample}).

With more work, in Theorem~\ref{thm:mainp-ample}, we generalize
to the case where $\calH$ is replaced with an $F$-ample coherent sheaf
tensored with a $p$-ample coherent sheaf. (See Definition~\ref{def:p-ample}
for the definition of $p$-ample.) This allows us to prove

\begin{theorem}\label{thm:equivalent}
Let $X$ be a projective scheme of pure dimension $d$, smooth over a field $k$. Let $\F_n$
be a sequence of coherent sheaves. Then the following are equivalent:
\begin{enumerate}
\item\label{equiv1} For any coherent $\G$, there exists $n_0$ such that $\G \otimes \F_n$
is generated by global sections for $n \geq n_0$.
\item\label{equiv2} For any coherent locally free $\calE$, there exists $n_1$ such that $\calE \otimes \F_n$
is generated by global sections for $n \geq n_1$.
\item\label{equiv3} For any invertible sheaf $\calH$, there exists $n_2$ such that
$\calH \otimes \F_n$ is $p$-ample for $n \geq n_2$.
\end{enumerate}
\end{theorem}

This answers \cite[Question~7.5]{Ke-filters-Frobenius}, at least when $X$
is smooth.

\section{Reductions}\label{sec:reductions}

In this section, we verify that Theorems~\ref{thm:FujitaConj}, \ref{thm:mainF-ample},
and \ref{thm:mainp-ample}
can be reduced to the case of $k$ algebraically closed and 
$\character k = p > 0$. 
First, we will define $F$-ample and $p$-ample.
Beyond these definitions, the material is standard, 
yet somewhat scattered in the literature.
Thus the familiar reader may wish to move to the next section. 

We recall the definition of $F$-ample  \cite{AraK} (also known as ``cohomologically $p$-ample''
\cite{Hart-ample}).

\begin{definition}\label{def:F-ample}
Let $X$ be a projective scheme over a field $k$, and let $\F$ be a coherent sheaf.
If $\character k = p > 0$ and $F$ is the absolute Frobenius morphism, then
define $\F^{(p^n)} = F^{*n}\F$.
The sheaf $\F$ is \emph{$F$-ample} if 
for any locally free coherent sheaf $\calE$, there exists $n_0$ such that
\[
H^q(X, \calE \otimes \F^{(p^n)}) = 0,\quad q > 0, n \geq n_0.
\]
If $\character k = 0$, then $\F$ 
is $F$-ample if and only if $\F_s$ is $F$-ample
for all closed fibers on some arithmetic thickening.
\end{definition}

Note that if $\F$ is $F$-ample and $\LL$ is an ample line bundle, then
$\F \otimes \LL$ is $F$-ample \cite[Theorem~4.5]{AraK}.

Another possible definition for ``ampleness'' of a coherent sheaf
is $p$-ampleness. While previously only defined for vector bundles
in characteristic $p$,
the definition easily extends.

\begin{definition}\label{def:p-ample}
Let $X$ be a projective scheme over a field $k$, and let $\F$ be a coherent sheaf.
If $\character k = p > 0$ and $F$ is the absolute Frobenius morphism, then
define $\F^{(p^n)} = F^{*n}\F$.
The sheaf $\F$ is \emph{$p$-ample} if 
for any locally free coherent sheaf $\calE$, there exists $n_0$ such that
\[
\calE \otimes \F^{(p^n)}
\]
is generated by global sections for $n \geq n_0$.
If $\character k = 0$, then $\F$ 
is $p$-ample if and only if $\F_s$ is $p$-ample
for all closed fibers on some arithmetic thickening.
\end{definition}

We chose the definition above because it immediately follows that an $F$-ample
coherent sheaf is $p$-ample (see Remark~\ref{rem:pAmpleIsFAmple}). 
However, we will see that we could have allowed $\calE$
to be any coherent sheaf (see Remark~\ref{rem:pampleAnySheaf}). 

Also, note that, in general, $F$-ample is a stronger condition
than $p$-ample. Specifically, on $\PP^n, n\geq 2$, the tangent
bundle is $p$-ample but not $F$-ample \cite[Example~5.9]{AraK}.

We now move to proving our reductions to
the case of $k$ algebraically closed and 
$\character k = p > 0$. 
First, we must check that we can extend our field $k$.

\begin{lemma}\label{lem:faithfully-flat}
Let $X$ be a projective scheme of pure dimension $d$, smooth over
a field $k$. Let $k \to k'$ be a field extension. Let $X' = X\times_k k'$
and $\pi:X' \to X$
be the projection. Then
\begin{enumerate}
\item $X'$ is of pure dimension $d$ and smooth over $k'$,
\item\label{flat-omega} $\omega_{X'/k'} \cong \pi^*\omega_{X/k}$,
\item\label{flat-ample} a line bundle $\LL$ on $X$ is ample (resp. very ample)
if and only if $\pi^*\LL$ is ample (resp. very ample),
\item a coherent sheaf $\F$ is generated by global sections
(resp. $F$-ample, $p$-ample)
if and only if $\pi^*\F$ is generated by global sections
(resp. $F$-ample, $p$-ample).
\end{enumerate}
\end{lemma}
\begin{proof} 
Since $X$ is smooth over $k$ and of pure dimension $d$, we have
that $\Omega_{X/k}$ is locally free of rank $d$ \cite[VII.5.1]{AltK}.
Now $\pi^*\Omega_{X/k} \cong \Omega_{X'/k'}$ \cite[p.~110]{AltK}.
So $\Omega_{X'/k'}$ is also
locally free of rank $d$ and we must have that $X'$
is smooth and of pure dimension $d$ \cite[VII.5.3]{AltK}.
We also immediately have \eqref{flat-omega}.

The claim \eqref{flat-ample} is \cite[$\mathrm{IV}_2$, 2.7.2]{EGA}.

If $\oplus \OO_X \to \F$ is a surjection, then 
$\oplus \pi^*\OO_X \cong \oplus \OO_{X'} \to \pi^*\F$ is a surjection.
On the other hand, since $\pi:X' \to X$ is faithfully flat,
if $H^0(\F) \otimes_k \OO_X \to \F$ has non-zero cokernel,
then there is still a non-zero cokernel upon pulling back.
This proves the claim about globally generated sheaves.

If $\character k = 0$, then the claim regarding $F$-ampleness
or $p$-ampleness is clear by definition.

Consider the case of $\character k = p>0$.
The absolute Frobenius morphism commutes with any morphism
of ringed spaces. 
Thus, $\pi^*(\F^{(p^n)}) = (\pi^*\F)^{(p^n)}$.
So the claim regarding $F$-ampleness follows from
\cite[Lemma~3.8]{Ke-filters-Frobenius}.

To show a sheaf $\F$ is $p$-ample, it is sufficient
to show that for a fixed very ample sheaf $\OO_X(1)$, we have
that for any $b \in \NN$, $\OO_X(-b) \otimes \F^{(p^n)}$
is generated by global sections for $n \gg 0$.
But this happens if and only if $\pi^*\OO_X(-b) \otimes \pi^*\F^{(p^n)}$
is generated by global sections.
Given \eqref{flat-ample},
we have the claim about $p$-ampleness.
\end{proof}

We now consider the case of $\character k = 0$. Our main tool
is arithmetic thickening. That is, given $X \to \Spec k$, we
can find a sub-algebra $A$ of $k$, finitely generated over $\ZZ$, 
and a scheme $\tilde{X}$ with flat morphism to $\Spec A$ such that 
$\tilde{X} \times_A k = X$. By shrinking $A$ via localization, any 
finite diagram of coherent sheaves on $X$ can be translated to
a diagram of coherent sheaves on $\tilde{X}$. 
We refer the reader to \cite[$\S$1]{AraK} for more on thickenings.

Since $A$ is finitely generated over $\ZZ$, the residue fields
at the closed points of $\Spec A$ are finite fields, hence
perfect of characteristic $p > 0$. By working on closed fibers, results
on the generic fiber and hence $X$ can often be deduced.

\begin{lemma}
Let $X$ be a projective scheme of pure dimension $d$, 
smooth over a field $k$ with
$\character k = 0$. Let $\LL$ be an ample line bundle,
and let $\F$ be a globally generated coherent sheaf.
Then there exists an arithmetic thickening
$f:\tilde{X} \to S=\Spec A, \tilde{\LL}, \tilde{\F}$ such that
at every (closed) fiber $\tilde{X}_s$ of $f$
\begin{enumerate}
\item $\tilde{\LL}_s$ is  ample,
\item $\tilde{\F}_s$ is globally generated,
\item the morphism $f$ is smooth,
\item each $X_s$ is a smooth scheme of pure dimension $d$,
\item\label{omega1} $\omega_{\tilde{X}/A} \otimes_A k = \omega_X$,
\item\label{omega2} $\omega_{\tilde{X}/A} \otimes_A k(s) = \omega_{X_s}$.
\end{enumerate}
\end{lemma}
\begin{proof}
Let $k(\nu)$ be the residue
field at the generic point of $S$. Given Lemma~\ref{lem:faithfully-flat},
we may replace $k$ with $k(\nu)$ and $X$ with $X\times_A k(\nu)$.

Now
we can shrink $S$ so that $\tilde{\LL}$ is ample  
\cite[$\mathrm{{III}_1}$, 4.7.1]{EGA}.
And $S$
can be shrunk again so there is surjection 
$H^0(\tilde{\F}_s) \otimes_A \OO_{X_s} \to \tilde{\F}_s$ at each fiber $X_s$
\cite[$\mathrm{IV}_3$, 9.4.2]{EGA}.

We can  shrink $S$ to make $f$ smooth 
at every fiber \cite[$\mathrm{IV}_3$, 12.2.4]{EGA}.
Hence the morphism $f$ is smooth over $\Spec A$ \cite[VII.1.8]{AltK}.

Since $f$ is smooth, the sheaf of differentials $\Omega_{\tilde{X}/A}$
is locally free \cite[VII.5.1]{AltK}. 
Shrinking $S$ again, we can make the rank
of $\Omega_{\tilde{X}/A}$ constant \cite[$\mathrm{IV}_2$, 2.5.2]{EGA},
equal to $d = \dim X$.
This makes each fiber a smooth scheme of pure dimension $d$.
The claims \eqref{omega1} and \eqref{omega2}
are true for $\Omega_{\tilde{X}/A}$ by the discussion on \cite[p.~110]{AltK}.
Hence they are also true for the $\omega_{\tilde{X}/A} = \wedge^d \Omega_{\tilde{X}/A}$.
\end{proof}

Note that if a coherent sheaf $\F$ over a field
of characteristic $0$ is $F$-ample or $p$-ample, then, by definition,
$\F_s$ is 
$F$-ample or $p$-ample on each closed fiber.
Thus, all hypotheses of Theorems~\ref{thm:FujitaConj}, \ref{thm:mainF-ample}, and
\ref{thm:mainp-ample}
can be translated to the closed fibers of an arithmetic thickening.

Finally, we need that the conclusions of 
Theorems~\ref{thm:FujitaConj}, \ref{thm:mainF-ample}, and \ref{thm:mainp-ample}
on a closed fiber imply the conclusion for our original scheme
over a field of characteristic $0$.

\begin{lemma}\label{lem:bad-lemma}
\footnote{This lemma is incorrect. See the corrected Lemma~\ref{E-lem:very-ample}. }
\sout{
Let $X$ be a projective scheme  over a field $k$ with
$\character k = 0$. Let $\LL$ be a line bundle and $\F$ a coherent sheaf.
Then there exists an arithmetic thickening
$f:\tilde{X} \to S=\Spec A, \tilde{\LL},\tilde{\F}$. }
\begin{enumerate}
\item \sout{ If $\tilde{\LL}_s$ is very ample on some closed fiber $X_s$,
then $\LL$ is very ample. }
\item  \sout{ If $\tilde{\F}_s$ is globally generated on some closed fiber
$X_s$, then $\F$ is globally generated. }
\end{enumerate} 
\end{lemma}
\begin{proof}
Again we may replace $k$ with $k(\nu)$.
Since $S$ is affine, the discussion following 
\cite[$\mathrm{{III}_1}$, 4.7.1]{EGA} proves the statement regarding $\LL$.

If $\oplus \OO_{X_s} \to \F_s$ is surjective, then Nakayama's Lemma
shows that $\oplus \OO_{\tilde{X}} \to \tilde{\F}$ is surjective
at the stalk at $s$. But then the map is surjective at the generic point
as well.
\end{proof}

Note that since the definition of $p$-ample in characteristic $0$
depends on \emph{all} closed fibers for some thickening, we cannot reduce
to the characteristic $p$ case to prove Theorem~\ref{thm:equivalent}. However,
we will see that the proof is quite easy given Theorem~\ref{thm:mainp-ample}.

\section{Fujita's Conjecture for $F$-ample sheaves}

In this section, we shall prove Theorem~\ref{thm:mainF-ample}.

We shall use \emph{Castelnuovo-Mumford regularity} of a sheaf. 
Let $\LL$ be an ample, globally generated line bundle. Recall that
a coherent sheaf $\G$ is $m$-regular (with respect to $\LL$) if
\[
H^q(X, \G \otimes \LL^{m-q}) = 0, \quad q > 0.
\]
We only need the fact that if $\G$ is $m$-regular, then $\G \otimes \LL^m$
is generated by global sections \cite[1.8.5, 1.8.14]{PAG}. While
regularity is usually defined for schemes over an algebraically closed
field, the reductions of Lemma~\ref{lem:faithfully-flat}
show that we may work over any field. 

\begin{remark}\label{rem:pAmpleIsFAmple}
Castelnuovo-Mumford regularity easily shows that an $F$-ample coherent sheaf is
$p$-ample.
To see this, if $\F$ is $F$-ample, one has that
$\calE \otimes \F^{(p^n)}$ is $0$-regular, 
and hence globally generated,
for $n \gg 0$.
\end{remark}

We are interested in the regularity of direct images under powers of the
absolute Frobenius $F$.

\begin{lemma}\label{lem:FrobRegularity}
Let $X$ be a projective scheme over a perfect field $k$ 
of positive characteristic, and let $\LL$
be an ample, globally generated line bundle. For any locally free coherent sheaf $\calE$
and $F$-ample coherent sheaf $\F$, there exists $n_0$ such that
\[
F_*^n(\calE) \otimes \F
\]
is $\dim X$-regular (with respect to $\LL$) for $n \geq n_0$.
Hence $F_*^n(\calE) \otimes \F \otimes \LL^{\dim X}$ is generated
by global sections.
\end{lemma}
\begin{proof}
Let $d = \dim X$ and write $\LL = \OO(1)$ (though $\LL$
may not be very ample).
If $q > d$, then it is trivial that $H^q(F_*^n(\calE) \otimes \F(d-q))=0$.

Now for $q=1, \dots, d$, we have that $\F(d-q)$ is $F$-ample.
So there exists $n_0$ such that
$H^q(\calE \otimes \F^{(p^n)}(p^nd-p^nq)) = 0$ for $n \geq n_0$.

Now since $k$ is perfect, $F$ is finite, so
$H^q(F^n_* (\calE \otimes \F^{(p^n)}(p^nd-p^nq))) = 0, q > 0$
\cite[Exercise~III.4.1]{Hart}.
But since $F$ is finite, 
\[
F^n_* (\calE \otimes \F^{(p^n)}(p^nd-p^nq)) \cong 
F^n_* (\calE) \otimes \F(d-q)
\]
by \cite[Lemma~5.7]{AraK}. So we have
shown the $d$-regularity we desired.
\end{proof}

Finally, we need a certain presumably well-known exact sequence
of locally free sheaves. Assuming $X$ is integral, a proof was
presented in \cite[Lemma~4.4]{Ke-filters-Frobenius}, so
we will only indicate how to reduce to that proof.

\begin{lemma}\label{lem:flat-sequence}
Let $X$ be a projective scheme of pure dimension $d$, 
smooth over a perfect field of 
characteristic $p > 0$. Let $F:X \to X$ be the absolute Frobenius morphism.
For any $n > 0$, there is an exact sequence of coherent, locally free sheaves
\[
0 \to \calK \to F^n_*\omega_X \to \omega_X \to 0
\]
where $\omega_X = \wedge^d \Omega_{X/k}$ is the dualizing sheaf of $X$.
\end{lemma}
\begin{proof}
The map $ F^n_*\omega_X \to \omega_X$ is dual to $\OO_X \to F^n_*\OO_X$
via Grothendieck-Serre duality \cite[Exercises~III.6.10, 7.2]{Hart}.
It is obviously a local question whether the cokernel of $\OO_X \to F^n_*\OO_X$
is locally free, so we may check this on the connected components of $X$.
But since $X$ is smooth, the connected components are the reduced, irreducible
components \cite[Remark~III.7.9.1]{Hart}. Thus we have reduced to the integral scheme case,
and we may proceed as in \cite[Lemma~4.4]{Ke-filters-Frobenius}.
Finally, the sheaf $\omega_X = \wedge^d \Omega_{X/k}$ is still the dualizing sheaf
for a smooth scheme of pure dimension $d$
\cite[I.4.6]{AltK}.

As an interesting alternative when $k$ is algebraically closed,
consider any open affine integral $U \subset X$.
The injection $\OO_U \to F^n_*\OO_U$ is split \cite[Proposition~1.1.6]{FrobSplit}.
Thus the cokernel of $\OO_X \to F^n_*\OO_X$
is locally free. The proof then proceeds as before.
\end{proof}

We may now turn to our main theorem. The proof is now quite simple.

\begin{theorem}\label{thm:mainF-ample}
Let $X$ be a projective  scheme of pure dimension $d$, smooth over a field $k$. Let $\LL$ be
an ample, globally generated line bundle. Let $\F$ be an $F$-ample coherent
sheaf. Then
\[
\omega_X \otimes \LL^{d} \otimes \F
\]
is generated by global sections.
\end{theorem}
\begin{proof}
By the reductions of Section~\ref{sec:reductions}, we may assume that
$k$ is perfect and of positive characteristic. 
We will have to allow $X$ to be reducible, but $X$
will be smooth of pure dimension $d$.

By Lemma~\ref{lem:FrobRegularity}, there exists $n$ such that
$F^n_*(\omega_X) \otimes \LL^{d} \otimes \F$ is generated by global sections.
But quotients of globally generated sheaves are globally generated, so
tensoring the exact sequence of
 Lemma~\ref{lem:flat-sequence} by $\LL^{d} \otimes \F$, 
 we have the theorem.
\end{proof}

We may now immediately prove Theorem~\ref{thm:FujitaConj}
as a corollary.

\begin{proof}[\textbf{Proof of Theorem~\ref{thm:FujitaConj}}]
Since 
ample line bundles are exactly the $F$-ample line bundles \cite[Lemma~2.4]{AraK}, 
we have Theorem~\ref{thm:FujitaConj}\eqref{fuj1}
as an immediate corollary. We will need to use another technique of
regularity to prove Theorem~\ref{thm:FujitaConj}\eqref{fuj2}.

By the reductions of Section~\ref{sec:reductions}, we may assume that
$k$ is algebraically closed and of positive characteristic.
Let $\calN = \omega_X \otimes \LL^{\dim X} \otimes \calH$,
which is globally generated by Theorem~\ref{thm:FujitaConj}(\ref{fuj1}).

Since $k$ is algebraically closed
and $\calN \otimes \LL$ is globally generated, it is sufficient to show
that $\scrm_x \otimes \calN \otimes \LL$
is generated by global sections for any closed $x \in X$
\cite[p.~21]{Hart-amplesub}. To show this, we follow the method of
\cite[Example~1.8.22]{PAG}.

Fix a closed point $x \in X$. Since $\LL$ is ample
and generated by global sections, we can let $Z$ be a zero-dimensional scheme
cut out by $d$ general sections of $\LL$ that vanish at $x$,
where $d = \dim X$.
Let $\calI_Z$ be the ideal sheaf of $Z$
and let $\calE = \oplus_{i=1}^{d} \LL$. Now $\calI_Z$ has a resolution
derived from the Koszul complex of $Z$ \cite[Example~1.8.18, Appendix~B.2]{PAG}:
\[
0 \to \wedge^d \calE^* \to \dots \wedge^2 \calE^* \to \calE^* \to \calI_Z \to 0.
\]
The exterior power $\wedge^q \calE^*$ is a direct sum of $\LL^{-q}$.
By Lemma~\ref{lem:FrobRegularity}, there exists $n$ such that 
$F_*^n(\omega_X) \otimes \calH$ is $\dim X$-regular.
So we can tensor the Koszul complex
with the locally free sheaf $F^n_*(\omega_X)
\otimes \LL^{\dim X + 1} \otimes \calH$, chase through
long exact sequences, and find that 
$\calI_Z \otimes F^n_*(\omega_X)
\otimes \LL^{\dim X + 1} \otimes \calH$ is $0$-regular.

Any sheaf with $0$-dimensional support is trivially $0$-regular.
So from the exact sequence
\[
0 \to \calI_Z \to \scrm_x \to \scrm_x/\calI_Z \to 0
\]
we get that $\scrm_x \otimes F^n_*(\omega_X)
\otimes \LL^{\dim X + 1} \otimes \calH$ is also $0$-regular
and thus generated by global sections.
The quotient sheaf $\scrm_x \otimes \calN \otimes \LL$ is then
generated by global sections and hence $\calN \otimes \LL$ is very ample.
\end{proof}

\begin{remark}
The above proof actually yields a  statement stronger than
Theorem~\ref{thm:FujitaConj}\eqref{fuj2}. Let $x \in X$.
Then
\[
\scrm_x \otimes \LL^{d+1} \otimes \F
\]
is generated by global sections for any locally free, $F$-ample
sheaf $\F$.
\end{remark}

\section{Fujita's Conjecture for $p$-ample sheaves}

We now can prove another variant of Fujita's Conjecture
regarding $p$-ample sheaves.

\begin{theorem}\label{thm:mainp-ample}
Let $X$ be a projective scheme of pure dimension $d$, smooth over a field $k$. Let $\LL$ be
an ample, globally generated line bundle. Let $\F$ be an $F$-ample coherent
sheaf and let $\calP$ be a $p$-ample coherent sheaf. Then
\[
\omega_X \otimes \LL^{d} \otimes \F \otimes \calP
\]
is generated by global sections. In particular, $\omega_X \otimes \LL^{d+1}\otimes \calP$
is generated by global sections.
\end{theorem}
\begin{proof}
By the reductions of Section~\ref{sec:reductions}, we may assume that
$k$ is perfect and of positive characteristic. 

Since $\calP$ is $p$-ample, there exits $n_0$ such that
$\omega_X \otimes \calP^{(p^n)}$ is generated by global sections
for $n \geq n_0$. We can apply the exact functor $F_*^n$
to the surjection $\oplus \OO_X \to \omega_X \otimes \calP^{(p^n)}$.
By the 
projection formula for finite morphisms \cite[Lemma~5.7]{AraK},
there is a surjective homomorphism
\begin{equation}\label{surj}
\oplus F_*^n(\OO_X) \to F_*^n(\omega_X) \otimes \calP \to 0
\end{equation}
for $n \geq n_0$.

Now by Lemma~\ref{lem:FrobRegularity}, $F_*^n(\OO_X) \otimes \LL^{\dim X} \otimes \F$
is $0$-regular with respect to $\LL$, and hence is generated by global sections,
for $n \gg 0$.
From Lemma~\ref{lem:flat-sequence} 
and the surjection \eqref{surj}, the quotient sheaf 
$\omega_X \otimes \LL^{\dim X} \otimes \F \otimes \calP$ is also
generated by global sections.
\end{proof}

As noted in the introduction, Theorem~\ref{thm:FujitaConj} was already
known for $k=\CC$, via the Kodaira Vanishing Theorem and
Castelnuovo-Mumford regularity. There is a Kodaira Vanishing Theorem
for $F$-ample sheaves (that is, for
$F$-ample $\F$,
$H^q(\omega_X \otimes \F) = 0, q>0$ \cite[Corollary~8.6]{AraK}).
Thus in characteristic $0$, Theorem~\ref{thm:mainF-ample}
could be proven via this ``Kodaira Vanishing--regularity'' method. 
However, \cite[Remark~7.4]{Ke-filters-Frobenius}
implies that there cannot be a Kodaira Vanishing Theorem for
$p$-ample sheaves, so that method is not an option for Theorem~\ref{thm:mainp-ample}.

We now turn to Theorem~\ref{thm:equivalent}. Let $\{ \F_n \}$ be
a sequence of coherent sheaves.
In \cite{Ke-filters-Frobenius},
the author thoroughly examined under what conditions $\calE \otimes \F_n$
had vanishing higher cohomology for $n \gg 0$. It was found that this
vanishing occurs if and only if for any invertible sheaf $\calH$, we
have that $\calH \otimes \F_n$ is $F$-ample for $n \gg 0$ 
\cite[Theorem~1.3]{Ke-filters-Frobenius}.
We now prove the global generation analogue.

\begin{proof}[\textbf{Proof of Theorem~\ref{thm:equivalent}}]
That \eqref{equiv1} implies \eqref{equiv2} is trivial.
Assuming \eqref{equiv2}, let $\LL$ be an ample invertible sheaf
and let $\calH$ be any invertible sheaf.
For $n \gg 0$, we have that $\LL^{-1} \otimes \calH \otimes \F_n$ is generated by global
sections. Thus $\calH \otimes \F_n$ is a quotient of $\oplus \LL$. It is easy to
see that direct sums and quotients of $p$-ample sheaves are $p$-ample.
Hence we have \eqref{equiv3}.

Now assume \eqref{equiv3}. Let $\G$ be a coherent sheaf and let $\LL$
be an invertible sheaf, ample and generated by global sections. There exits $m$
such that $\G \otimes \LL^m$ is generated by global sections.
By the assumption on $\{ \F_n \}$, there exists $n_0$ such that
$\omega_X^{-1} \otimes \LL^{-m - \dim X - 1} \otimes \F_n$ is $p$-ample
for $n \geq n_0$. Then by Theorem~\ref{thm:mainp-ample}, \[
\LL^{-m} \otimes \F_n \cong (\omega_X \otimes \LL^{\dim X + 1})
\otimes (\omega_X^{-1} \otimes \LL^{-m - \dim X - 1} \otimes \F_n)
\]
is generated by global sections for $n \geq n_0$. 
Thus $\G \otimes \F_n$ is generated by global sections,
as it is a tensor product of globally generated sheaves.
\end{proof}

\begin{remark}\label{rem:pampleAnySheaf}
The fact that \ref{thm:equivalent}\eqref{equiv2} implies \ref{thm:equivalent}\eqref{equiv1}
shows that we could have allowed the $\calE$ in Definition~\ref{def:p-ample}
to be any coherent sheaf.
On the other hand, the analogue for $F$-ample sheaves is not known to
be true, unless the $F$-ample sheaf only fails to be locally free
on a subscheme of dimension $\leq 2$ \cite[Lemma~3.11]{Ke-filters-Frobenius}.
\end{remark}

\begin{remark}
Instead of just taking a sequence of coherent sheaves, we could
have indexed by a \emph{filter} (that is, a partially ordered set
such that for any $\alpha, \beta$, there exist $\gamma$ with $\alpha < \gamma,
\beta < \gamma$).
Thus, Theorem~\ref{thm:equivalent} answers \cite[Question~7.5]{Ke-filters-Frobenius}
in the affirmative, at least when $X$ is smooth. Based on the results
of \cite[Theorem~7.2]{Ke-filters-Frobenius}, we conjecture that Theorem~\ref{thm:equivalent}
will remain true for any projective $X$, at least when the $\F_n$
are locally free.  However, we do not see a way to reduce to the smooth case,
as in \cite[Section~5]{Ke-filters-Frobenius}.
\end{remark}

\noindent
\textbf{Acknowledgement}. We would like to thank the referee for helpful comments.

\bibliographystyle{amsplain}
\providecommand{\bysame}{\leavevmode\hbox to3em{\hrulefill}\thinspace}

\newpage

\numberwithin{equation}{section}
\renewcommand{\theequation}{E\arabic{section}.\arabic{equation}}
\setcounter{section}{0}
\renewcommand*{\theHsection}{E.\the\value{section}}

\begin{center}
\LARGE Erratum
\end{center}

\section{Corrected Lemma}
\label{SecErratum}

Consider an arithmetic thickening 
$f:\tilde{X} \to S=\Spec A$.
In \cite[Lemma~\ref{lem:bad-lemma}]{E-Ke}, the author incorrectly treated 
very ample and globally generated as properties
preserved on fibers $\tilde{X}_s$ where $s$ ranges over an \emph{open} subset $T$ of the base
scheme $S$.
However, such sets $T$ are only \emph{constructible}. Indeed, if globally generated was
an open condition on the Noetherian base $S$, then the same would hold for semi-ample
(that is, some power of a line bundle would be globally generated).
However, \cite[Theorem~3.0]{E-Keel} gave an example
of a naturally defined line bundle $\LL$ which is semi-ample in positive characteristic, but
not in characteristic zero. Hence, $\tilde{\LL}$ would be semi-ample
at closed points of an arithmetic thickening, but not at the generic point.

Fortunately, the results proven in \cite{E-Ke} hold over all closed points
simultaneously. Thus this corrected lemma is a satisfactory replacement
for \cite[Lemma~2.5]{E-Ke}.

\begin{lemma}\label{E-lem:very-ample}
Let $X$ be a projective scheme  over a field $k$ with
$\character k = 0$. Let $\LL$ be a line bundle and $\F$ a coherent sheaf on $X$.
Let
$f:\tilde{X} \to S=\Spec A, \tilde{\LL},\tilde{\F}$ be an arithmetic thickening.
\begin{enumerate}
\item\label{part1} If $\tilde{\LL}_s$ is very ample on every closed fiber $X_s$,
then $\LL$ is very ample.
\item\label{part2} If $\tilde{\F}_s$ is globally generated on every closed fiber
$X_s$, then $\F$ is globally generated.
\end{enumerate}
\end{lemma}
\begin{proof}
By definition of arithmetic thickening, the ring $A$ is a finite-type
$\ZZ$-algebra and hence a Jacobson ring 
\cite[Tag 00GC]{E-stacks-project}.
So $S$ is a Jacobson space \cite[Tag 00G3]{E-stacks-project}. 
Any open subset $U$ is also a Jacobson space and the closed
points of $U$ are closed in $S$ \cite[Tag 005X]{E-stacks-project}.
Therefore, we may replace $A$ with any one element localization
and still keep the given hypotheses on all $\tilde{\LL}_s$
or $\tilde{\F}_s$.

Our goal is to show (possibly after shrinking $S$)
 that if $T$ is the set of (not necessarily closed) points $s$
where $\tilde{\LL}_s$ is very ample or where $\tilde{\F}_s$ is globally generated
on $X_s$, then $T$ is a constructible set. Assume this is done. The set of all closed points
of a Jacobson space is a dense set by definition \cite[Tag 005U]{E-stacks-project},
so $T$ is a dense constructible subset of $S$. 
Since $A$ is an integral domain, the space $S$ is irreducible
and so $T$ contains the generic point $\eta$ of $S$
\cite[Tag 005K]{E-stacks-project}. Then we are done by 
\cite[Lemma~2.3]{E-Ke}, noting that the smooth hypothesis was not needed
for parts 3 and 4 of that lemma.

We begin by shrinking $S$ until $f$ is projective \cite[$\mathrm{IV}_3$, 8.10.5]{E-EGA}.
If $T$ is the set of points $s \in S$
such that $\tilde{\LL}_s$ is very ample, then $T$ is
constructible \cite[$\mathrm{IV}_3$, 9.6.2]{E-EGA}. This proves part \eqref{part1}.

Now by Generic Flatness, we may assume $f$ is flat and $\tilde{\F}$ is flat over $S$
\cite[Tag 052A]{E-stacks-project}. Further, we may assume $s \to h^0({X}_s, \tilde{\F}_s)$
is constant on $S$ \cite[Tag 0BDN]{E-stacks-project} and so $H^0(\tilde{X}, \tilde{\F})$
is a locally free $A$-module and
$H^0(\tilde{X}, \tilde{\F}) \otimes k(s) \to H^0({X}_s, \tilde{\F}_s)$
is an isomorphism for all $s \in S$ \cite[III.12.9]{E-Hartshorne}.
We can shrink $S$ further to assume $H^0(\tilde{X}, \tilde{\F})$
is a free $A$-module of rank $r = h^0({X}_\eta, \tilde{\F}_\eta)$.

Consider the natural morphism of $\OO_{\tilde{X}}$-modules 
$u:f^* f_* \tilde{\F} \to \tilde{\F}$. Note that
$f_* \tilde{\F} \cong \OO_S^r$ 
and so pulling back to a fiber ${X}_s$, we have
\[
u_s: H^0({X}_s, \tilde{\F}_s) \otimes \OO_{{X}_s} \to \tilde{\F}_s.
\]
By hypothesis, $u_s$ is surjective for all closed points $s \in S$. 
But the set of all points with $u_s$ surjective
is constructible \cite[$\mathrm{IV}_3$, 9.4.5]{E-EGA}.
Thus, $\eta$ is an element of this constructible set and we have \eqref{part2}.
\end{proof}

\end{document}